\documentclass[10pt]{amsart}
\usepackage{amsfonts}
\usepackage{amsmath}
\usepackage{amsthm}
\usepackage{amssymb}
\usepackage{latexsym}
\newtheorem{cor}{Corollary}[section]
\newtheorem{lem}{Lemma}[section]

\newtheorem{prop}{Proposition}[section]

\theoremstyle{definition}
\newtheorem{defn}{Definition}[section]
\theoremstyle{definition}
\newtheorem{thm}{Theorem}

\newtheorem*{rem}{Remark}

\newenvironment{pf}{\proof}{\endproof}

\theoremstyle{remark}

\numberwithin{equation}{section}
\begin{document}

\newcommand{\thmref}[1]{Theorem~\ref{#1}}
\newcommand{\secref}[1]{Sect.~\ref{#1}}
\newcommand{\lemref}[1]{Lemma~\ref{#1}}
\newcommand{\propref}[1]{Proposition~\ref{#1}}
\newcommand{\corref}[1]{Corollary~\ref{#1}}
\newcommand{\remref}[1]{Remark~\ref{#1}}
\newcommand{\nc}{\newcommand}
\newcommand{\rnc}{\renewcommand}
\nc{\cal}{\mathcal}
\nc{\goth}{\mathfrak}
\rnc{\bold}{\mathbf}
\renewcommand{\frak}{\mathfrak}
\renewcommand{\Bbb}{\mathbb}

\nc{\Cal}{\mathcal}
\nc{\Xp}[1]{X^+(#1)}
\nc{\Xm}[1]{X^-(#1)}
\nc{\on}{\operatorname}
\nc{\ch}{\mbox{ch}}
\nc{\Z}{{\bold Z}}
\nc{\J}{{\cal J}}
\nc{\C}{{\bold C}}
\nc{\Q}{{\bold Q}}
\renewcommand{\P}{{\cal P}}
\nc{\N}{{\Bbb N}}
\nc\beq{\begin{equation}}
\nc\enq{\end{equation}}
\nc\lan{\langle}
\nc\ran{\rangle}
\nc\bsl{\backslash}
\nc\mto{\mapsto}
\nc\lra{\leftrightarrow}
\nc\hra{\hookrightarrow}
\nc\sm{\smallmatrix}
\nc\esm{\endsmallmatrix}
\nc\sub{\subset}
\nc\ti{\tilde}
\nc\nl{\newline}
\nc\fra{\frac}
\nc\und{\underline}
\nc\ov{\overline}
\nc\ot{\otimes}
\nc\bbq{\bar{\bq}_l}
\nc\bcc{\thickfracwithdelims[]\thickness0}
\nc\ad{\text{\rm ad}}
\nc\Ad{\text{\rm Ad}}
\nc\Hom{\text{\rm Hom}}
\nc\End{\text{\rm End}}
\nc\Ind{\text{\rm Ind}}
\nc\Res{\text{\rm Res}}
\nc\Ker{\text{\rm Ker}}
\rnc\Im{\text{Im}}
\nc\sgn{\text{\rm sgn}}
\nc\tr{\text{\rm tr}}
\nc\Tr{\text{\rm Tr}}
\nc\supp{\text{\rm supp}}
\nc\card{\text{\rm card}}
\nc\bst{{}^\bigstar\!}
\nc\he{\heartsuit}
\nc\clu{\clubsuit}
\nc\spa{\spadesuit}
\nc\di{\diamond}

\nc\al{\alpha}
\nc\bet{\beta}
\nc\ga{\gamma}
\nc\de{\delta}
\nc\ep{\epsilon}
\nc\io{\iota}
\nc\om{\omega}
\nc\si{\sigma}
\rnc\th{\theta}
\nc\ka{\kappa}
\nc\la{\lambda}
\nc\ze{\zeta}

\nc\vp{\varpi}
\nc\vt{\vartheta}
\nc\vr{\varrho}

\nc\Ga{\Gamma}
\nc\De{\Delta}
\nc\Om{\Omega}
\nc\Si{\Sigma}
\nc\Th{\Theta}
\nc\La{\Lambda}
\nc\boa{\bold a}
\nc\bob{\bold b}
\nc\boc{\bold c}
\nc\bod{\bold d}
\nc\boe{\bold e}
\nc\bof{\bold f}
\nc\bog{\bold g}
\nc\boh{\bold h}
\nc\boi{\bold i}
\nc\boj{\bold j}
\nc\bok{\bold k}
\nc\bol{\bold l}
\nc\bom{\bold m}
\nc\bon{\bold n}
\nc\boo{\bold o}
\nc\bop{\bold p}
\nc\boq{\bold q}
\nc\bor{\bold r}
\nc\bos{\bold s}
\nc\bou{\bold u}
\nc\bov{\bold v}
\nc\bow{\bold w}
\nc\boz{\bold z}

\nc\ba{\bold A}
\nc\bb{\bold B}
\nc\bc{\bold C}
\nc\bd{\bold D}
\nc\be{\bold E}
\nc\bg{\bold G}
\nc\bh{\bold h}
\nc\bH{\bold H}

\nc\bi{\bold I}
\nc\bj{\bold J}
\nc\bk{\bold K}
\nc\bl{\bold L}
\nc\bm{\bold M}
\nc\bn{\bold N}
\nc\bo{\bold O}
\nc\bp{\bold P}
\nc\bq{\bold Q}
\nc\br{\bold R}
\nc\bs{\bold S}
\nc\bt{\bold T}
\nc\bu{\bold U}
\nc\bv{\bold v}
\nc\bV{\bold V}

\nc\bw{\bold W}
\nc\bz{\bold Z}
\nc\bx{\bold x}
\nc\bX{\bold X}
\nc\blambda{{\mbox{\boldmath $\Lambda$}}}
\nc\bpi{{\mbox{\boldmath $\pi$}}}

\nc\e[1]{E_{#1}}
\nc\ei[1]{E_{\delta - \alpha_{#1}}}
\nc\esi[1]{E_{s \delta - \alpha_{#1}}}
\nc\eri[1]{E_{r \delta - \alpha_{#1}}}
\nc\ed[2][]{E_{#1 \delta,#2}}
\nc\ekd[1]{E_{k \delta,#1}}
\nc\emd[1]{E_{m \delta,#1}}
\nc\erd[1]{E_{r \delta,#1}}

\nc\ef[1]{F_{#1}}
\nc\efi[1]{F_{\delta - \alpha_{#1}}}
\nc\efsi[1]{F_{s \delta - \alpha_{#1}}}
\nc\efri[1]{F_{r \delta - \alpha_{#1}}}
\nc\efd[2][]{F_{#1 \delta,#2}}
\nc\efkd[1]{F_{k \delta,#1}}
\nc\efmd[1]{F_{m \delta,#1}}
\nc\efrd[1]{F_{r \delta,#1}}
\nc{\ug}{\bu^{fin}}

\nc\fa{\frak a}
\nc\fb{\frak b}
\nc\fc{\frak c}
\nc\fd{\frak d}
\nc\fe{\frak e}
\nc\ff{\frak f}
\nc\fg{\frak g}
\nc\fh{\frak h}
\nc\fj{\frak j}
\nc\fk{\frak k}
\nc\fl{\frak l}
\nc\fm{\frak m}
\nc\fn{\frak n}
\nc\fo{\frak o}
\nc\fp{\frak p}
\nc\fq{\frak q}
\nc\fr{\frak r}
\nc\fs{\frak s}
\nc\ft{\frak t}
\nc\fu{\frak u}
\nc\fv{\frak v}
\nc\fz{\frak z}
\nc\fx{\frak x}
\nc\fy{\frak y}

\nc\fA{\frak A}
\nc\fB{\frak B}
\nc\fC{\frak C}
\nc\fD{\frak D}
\nc\fE{\frak E}
\nc\fF{\frak F}
\nc\fG{\frak G}
\nc\fH{\frak H}
\nc\fJ{\frak J}
\nc\fK{\frak K}
\nc\fL{\frak L}
\nc\fM{\frak M}
\nc\fN{\frak N}
\nc\fO{\frak O}
\nc\fP{\frak P}
\nc\fQ{\frak Q}
\nc\fR{\frak R}
\nc\fS{\frak S}
\nc\fT{\frak T}
\nc\fU{\frak U}
\nc\fV{\frak V}
\nc\fZ{\frak Z}
\nc\fX{\frak X}
\nc\fY{\frak Y}
\nc\tfi{\ti{\Phi}}
\nc\bF{\bold F}

\nc\ua{\bold U_\A}

\nc\qinti[1]{[#1]_i}
\nc\q[1]{[#1]_q}
\nc\xpm[2]{E_{#2 \delta \pm \alpha_#1}}  
\nc\xmp[2]{E_{#2 \delta \mp \alpha_#1}}
\nc\xp[2]{E_{#2 \delta + \alpha_{#1}}}
\nc\xm[2]{E_{#2 \delta - \alpha_{#1}}}
\nc\hik{\ed{k}{i}}
\nc\hjl{\ed{l}{j}}
\nc\qcoeff[3]{\left[ \begin{smallmatrix} {#1}& \\ {#2}& \end{smallmatrix}
\negthickspace \right]_{#3}}
\nc\qi{q}
\nc\qj{q}

\nc\ufdm{{_\ca\bu}_{\rm fd}^{\le 0}}


\nc\isom{\cong} 

\nc{\pone}{{\Bbb C}{\Bbb P}^1}
\nc{\pa}{\partial}
\def\H{\cal H}
\def\L{\cal L}
\nc{\F}{{\cal F}}
\nc{\Sym}{{\goth S}}
\nc{\A}{{\cal A}}
\nc{\arr}{\rightarrow}
\nc{\larr}{\longrightarrow}

\nc{\ri}{\rangle}
\nc{\lef}{\langle}
\nc{\W}{{\cal W}}
\nc{\uqatwoatone}{{U_{q,1}}(\su)}
\nc{\uqtwo}{U_q(\goth{sl}_2)}
\nc{\dij}{\delta_{ij}}
\nc{\divei}{E_{\alpha_i}^{(n)}}
\nc{\divfi}{F_{\alpha_i}^{(n)}}
\nc{\Lzero}{\Lambda_0}
\nc{\Lone}{\Lambda_1}
\nc{\ve}{\varepsilon}
\nc{\phioneminusi}{\Phi^{(1-i,i)}}
\nc{\phioneminusistar}{\Phi^{* (1-i,i)}}
\nc{\phii}{\Phi^{(i,1-i)}}
\nc{\Li}{\Lambda_i}
\nc{\Loneminusi}{\Lambda_{1-i}}
\nc{\vtimesz}{v_\ve \otimes z^m}

\nc{\asltwo}{\widehat{\goth{sl}_2}}
\nc\eh{\frak h^e}  
\nc\loopg{L(\frak g)}  
\nc\eloopg{L^e(\frak g)} 
\nc\ebu{\bu^e} 
\nc\loopa{L(\frak a)}  

\nc\teb{\tilde E_\boc}
\nc\tebp{\tilde E_{\boc'}}

\title{Tensor products of level zero representations}\author{Vyjayanthi Chari}
\address{Vyjayanthi Chari, Department of Mathematics, University of California, 
Riverside, CA 92521.}
\maketitle
\section{Introduction} In this paper we  give a sufficient condition for the 
tensor product of irreducible finite--dimensional representations of quantum 
affine algebras to be cyclic. In particular, this proves a generalization of a 
recent result of Kashiwara \cite{K},   and also establishes a conjecture stated 
in that paper. 

We describe our results in some detail. Let ${\frak{g}}$ be a    complex simple 
finite--dimensional Lie algebra of rank $n$, and let $\bu_q$ be the quantized 
untwisted affine algebra over $\bc(q)$ associated to $\frak g$. For every  
$n$--tuple $\bpi=(\pi_1,\cdots ,\pi_n)$ of polynomials  with coefficients in 
$\bc(q)[u]$ and with constant term one, there exists  a unique (up to 
isomorphism)  irreducible finite--dimensional representation  $V(\bpi)$ of 
$\bu_q$. For each element  $w$ in the Weyl group $W$ of $\frak g$,  let 
$v_{w\bpi}$ be the extremal vector defined in \cite{K}.  
 In this paper we compute the action of the imaginary root vectors in $\bu_q$ on 
the elements $v_{w\bpi}$. To do this we define in Section 2 an action of the 
braid group $\cal{B}$  of $\frak g$ on elements of $(\bc(q)[[u]])^n$ and prove 
that the eigenvalue of $v_{w\bpi}$ is the element $T_w(\bpi)$ where 
$T:W\to\cal{B}$ is the canonical section defined in \cite{Bo}.

To state our result we assume  for simplicity (in the introduction only)  that 
$\frak g$ is simply--laced.   We shall again for simplicity, only deal with 
polynomials in $\bc(q)[u]$ which split into linear factors. Any such polynomial 
can be written uniquely as a product 
\begin{equation*}\pi(u)=\prod_{r=1}^k(1-a_rq^{m_r-1}u)(1-a_rq^{m_r-3}u)\cdots 
(1-a_rq^{-m_r+1}u),\end{equation*}
where $a_r\in\bc(q)$ and $m_r\in\bz_+$ satisfy
\begin{equation*}  \frac{a_r}{a_l}\ne q^{\pm (m_r+ m_l-2m)},\ \ 0\le 
m<\text{min}(m_r,m_l), \ \ \end{equation*}  if $r<l$. Let $S(\pi)$ be the 
collection of the pairs $(a_r,m_r)$. $1\le r\le k$ defined above.
  Say that a polynomial $\pi'(u) > \pi(u)$ if 
\begin{equation*}\frac{a_r'}{a_\ell}\ne q^{m_{r}'-m_\ell-2p}, \ \ \forall \ \ 1 
\le p\le m_r', \end{equation*}
for all pairs $(a_r',m_r')\in S(\pi')$ and $(a_\ell, m_\ell)\in S(\pi)$. 
Let $s_1, s_2,
\cdots ,s_n$ be the set of simple reflections in $W$. Our main result says that: 

\medspace 

\noindent {\it{ The tensor product $V(\bpi')\otimes V(\bpi)$ is cyclic on 
$v_{\bpi'}\otimes v_{\bpi}$ if, for all $w\in W$ and for all $i=1,\cdots ,n $ 
with $\ell(s_iw)=\ell(w)+1$, we have
\begin{equation*} (T_w\bpi')_i\ > \ \pi_i.\end{equation*} More generally, if 
$V_1,\cdots V_r$ are irreducible finite--dimensional representations, then 
$V_1\otimes\cdots \otimes V_r$ is cyclic, if every pair $V_j\otimes V_l$ is 
cyclic for all $j<\ell $}}

To make the connection with Kashiwara's theorem and conjectures, we consider the 
 case 
\begin{equation*} \pi_j(u) =1\ \ (j\ne i), \ \ 
\pi_i(u)=\prod_{s=1}^m(1-q^{m+1-2s}au)\ \ (a\in\bc(q)). \end{equation*}
Denoting this $n$--tuple of polynomials as $\bpi^i_{m,a}$ and the corresponding 
representation by $V_i(m,a)$, we prove the following:
\medspace
\noindent {\it{ Let $l\ge 1$ and let $i_j\in I$, $m_j\in\bz_+$, $a_j\in\bc(q)$  
for $1\le j\le l$. The tensor product $V_{i_1}(m_1,a_1)\otimes V_{i_2}(m_2, 
a_2)\otimes\cdots\otimes V_{i_l}(m_l,a_l)$ is cyclic on the tensor product of 
highest weight vectors if for all $r<s$, 
\begin{equation*} \frac{a_r}{a_s}\ne q^{m_r-m_s-p}, \ \ \forall \ p\ge 
0.\end{equation*}}}

\noindent The case when $m_j=1$ was originally conjectured  and partially proved 
in \cite{AK} and  completely proved in \cite{K} (and in \cite{VV} for the 
simply-laced case). The result in the case  when the $m_i$ are arbitrary but   
$a_i=1$  for all $i$ was conjectured in \cite{K}, \cite{HKOTY}.

\section{Preliminaries} In this section we recall the definition of 
quantum affine algebras  and several results on the classification of their 
irreducible  finite--dimensional representations.
 
Let $q$ be an indeterminate, let $\bc(q)$ be the field of rational
functions in $q$ with complex coefficients.  For $r,m\in\bn$, $m\ge r$, define
\begin{equation*} 
[m]_q=\frac{q^m -q^{-m}}{q -q^{-1}},\ \ \ \ [m]_q! =[m]_q[m-1]_q\ldots 
[2]_q[1]_q,\ \ \ \ 
\left[\begin{matrix} m\\ r\end{matrix}\right]_q 
= \frac{[m]_q!}{[r]_q![m-r]_q!}.
\end{equation*}

 Let $\frak g$ be a complex finite--dimensional simple Lie algebra  of rank $n$, 
 let $I=\{1,2,\cdots ,n\}$, let $\{\alpha_i:i\in I\}$ be the set of simple roots 
and let $\{\omega_i:i\in I\}$ be the set of fundamental weights. Let $Q^+$ 
(resp. $P^+$) be the non--negative root (resp. weight) lattice of $\frak g$. 
 Let $A=(a_{ij})_{i,j\in I}$  be the $n\times n$ Cartan matrix of $\frak g$ and  
let $\hat A 
=(a_{ij})$ be  the $(n+1)\times (n+1)$ extended Cartan matrix associated to 
$\frak g$. Let $\hat{I} =I\cup\{0\}$. Fix non--negative integers $d_i$  
($i\in\hat{I}$) such that the matix $(d_ia_{ij})$ is symmetric. Set 
$q_i=q^{d_i}$ and  $[m]_i=[m]_{q_i}$.
\begin{prop}{\label{defnbu}} There is a Hopf algebra $\tilde{\bu}_q$ over 
$\bc(q)$ which is generated as an algebra by elements $E_{\alpha_i}$, 
$F_{\alpha_i}$, $K_i^{{}\pm 1}$ ($i\in\hat I$), with the following defining 
relations:
\begin{align*} 
  K_iK_i^{-1}=K_i^{-1}K_i&=1,\ \ \ \ K_iK_j=K_jK_i,\\ 
  K_iE_{\alpha_j} K_i^{-1}&=q_i^{ a_{ij}}E_{\alpha_j},\\ 
K_iF_{\alpha_j} K_i^{-1}&=q_i^{-a_{ij}}F_{\alpha_j},\\
  [E_{\alpha_i}, F_{\alpha_j}
]&=\delta_{ij}\frac{K_i-K_i^{-1}}{q_i-q_i^{-1}},\\ 
  \sum_{r=0}^{1-a_{ij}}(-1)^r\left[\begin{matrix} 1-a_{ij}\\ 
  r\end{matrix}\right]_i
&(E_{\alpha_i})^rE_{\alpha_j}(E_{\alpha_i})^{1-a_{ij}-r}=0\ 
  \ \ \ \ \text{if $i\ne j$},\\
\sum_{r=0}^{1-a_{ij}}(-1)^r\left[\begin{matrix} 1-a_{ij}\\ 
  r\end{matrix}\right]_i
&(F_{\alpha_i})^rF_{\alpha_j}(F_{\alpha_i})^{1-a_{ij}-r}=0\ 
  \ \ \ \ \text{if $i\ne j$}.
\end{align*}
The comultiplication of $\tilde{\bu}_q$ is given on generators by
$$\Delta(E_{\alpha_i})=E_{\alpha_i}\ot 1+K_i\ot E_{\alpha_i},\ \ 
\Delta(F_{\alpha_i})=F_{\alpha_i}\ot K_i^{-1} + 1\ot F_{\alpha_i},\ \ 
\Delta(K_i)=K_i\ot K_i,$$
for $i\in\hat I$.\hfill\qedsymbol
\end{prop}

Set $K_{\theta} =\prod_{i=1}^n K_i^{r_i/d_i}$, where $\theta=\sum 
r_i\alpha_i$ is the highest root in $R^+$. 
Let $\bu_q$ be the  quotient of $\tilde{\bu}_q$ by the ideal generated by the 
central element  $K_0K_{\theta}^{-1}$; we call this the quantum loop algebra 
of $\frak g$.

It follows from 
\cite{Dr}, \cite{B}, \cite{J} that $\bu_q$ is isomorphic to the 
algebra with generators $x_{i,r}^{{}\pm{}}$ ($i\in I$, $r\in\bz$), $K_i^{{}\pm 
1}$ 
($i\in I$), $h_{i,r}$ ($i\in I$, $r\in \bz\backslash\{0\}$) and the following 
defining relations:
\begin{align*}
   K_iK_i^{-1} = K_i^{-1}K_i& =1, \ \  
 K_iK_j =K_jK_i,\\  K_ih_{j,r}& =h_{j,r}K_i,\\ 
 K_ix_{j,r}^\pm K_i^{-1} &= q_i^{{}\pm
    a_{ij}}x_{j,r}^{{}\pm{}},\ \ \\ 
  [h_{i,r},h_{j,s}]=0,\; \; & [h_{i,r} , x_{j,s}^{{}\pm{}}] =
  \pm\frac1r[ra_{ij}]_ix_{j,r+s}^{{}\pm{}},\\ 
 x_{i,r+1}^{{}\pm{}}x_{j,s}^{{}\pm{}} -q_i^{{}\pm
    a_{ij}}x_{j,s}^{{}\pm{}}x_{i,r+1}^{{}\pm{}} &=q_i^{{}\pm
    a_{ij}}x_{i,r}^{{}\pm{}}x_{j,s+1}^{{}\pm{}}
  -x_{j,s+1}^{{}\pm{}}x_{i,r}^{{}\pm{}},\\ [x_{i,r}^+ ,
  x_{j,s}^-]=\delta_{i,j} & \frac{ \psi_{i,r+s}^+ -
    \psi_{i,r+s}^-}{q_i - q_i^{-1}},\\ 
\sum_{\pi\in\Sigma_m}\sum_{k=0}^m(-1)^k\left[\begin{matrix}m\\k\end{matrix}
\right]_i
  x_{i, r_{\pi(1)}}^{{}\pm{}}\ldots x_{i,r_{\pi(k)}}^{{}\pm{}} &
  x_{j,s}^{{}\pm{}} x_{i, r_{\pi(k+1)}}^{{}\pm{}}\ldots
  x_{i,r_{\pi(m)}}^{{}\pm{}} =0,\ \ \text{if $i\ne j$},
\end{align*}
for all sequences of integers $r_1,\ldots, r_m$, where $m =1-a_{ij}$, $\Sigma_m$ 
is the symmetric group on $m$ letters, and the $\psi_{i,r}^{{}\pm{}}$ are 
determined by equating powers of $u$ in the formal power series 
$$\sum_{r=0}^{\infty}\psi_{i,\pm r}^{{}\pm{}}u^{{}\pm r} = K_i^{{}\pm 1} 
{\text{exp}}\left(\pm(q_i-q_i^{-1})\sum_{s=1}^{\infty}h_{i,\pm s} u^{{}\pm 
s}\right).$$ 
\vskip 12pt

For $i\in I$, the preceding  isomorphism maps $E_{\alpha_i}$ to  $x_{i,0}^+$ and 
$F_{\alpha_i}$ to $x_{i,0}^-$.  The subalgebra generated by $E_{\alpha_i}$, 
$F_{\alpha_i}$, $K_i^{\pm 1}$ ($i\in I$) is the quantized enveloping algebra 
$\bu_q^{fin}$ associated to $\frak g$. Let $\bu_q(<)$ be the subalgebra 
generated by the elements $x_{i,k}^-$ ($i\in I$, $k\in\bz\}$).  For $i\in I$, 
let $\bu_i$ be the subalgebra of $\bu_q$ generated by the elements 
$\{x_{i,k}^\pm: k\in\bz\}$, the subalgebra $\bu_i^{fin}$ is defined in the same 
way. Notice that $\bu_i$ is isomorphic  to  the quantum affine algebra 
$\bu_{q_i}(\widehat{sl_2})$.  Let $\Delta_i$ be the comultiplication of 
$\bu_{q_i}(\widehat{sl_2})$.

An explicit formula for the comultiplication on the Drinfeld generators is not 
known. Also, the subalgebra $\bu_i$ is not a Hopf subalgebra of $\bu_q$.  
However, partial information which is sufficient for our needs is  given in the 
next proposition.
Let \begin{equation*} X^\pm =\sum_{i\in I, k\in\bz}  \bc(q)x_{i,k}^\pm,\ \ 
X^\pm(i) =\sum_{j\in I\backslash\{i\}, k\in\bz}\bc(q) 
x_{j,k}^\pm.\end{equation*}

\begin{prop}\label{comultip} The restriction of  $\Delta$ to $\bu_i$ satisfies,
\begin{equation*} \Delta(x) =\Delta_i(x) \mod\left(\bu_q\otimes(\bu_q\backslash 
\bu_i)\right).\end{equation*}
More precisely:
\begin{enumerate}
\item[(i)] Modulo $\bu_qX^-\otimes \bu_q(X^+)^2 + \bu_q X^-\otimes \bu_qX^+(i)$, 
we have
\begin{align*} \Delta (x_{i,k}^+)& =x_{i,k}^+\otimes 1+ K_i\otimes x_{i,k}^+ 
+\sum_{j=1}^k \psi^+_{i,j}\otimes x_{i,k-j}^+ \ \ (k\ge 0),\\
\Delta (x_{i,-k}^+)& =K_i^{-1}\otimes x_{i,-k}^+ + x_{i,-k}^+\otimes 1 
+\sum_{j=1}^{k-1} \psi^-_{i, -j}\otimes  x_{i,-k+j}^+ \ \ (k> 0).\end{align*}
\item[(ii)]  Modulo $\bu_q(X^-)^2\otimes\bu_qX^+ + \bu_qX^-\otimes \bu_qX^+(i)$, 
we have
\begin{align*} \Delta (x_{i,k}^-)& =x_{i,k}^-\otimes K_i+ 1\otimes x_{i,k}^- 
+\sum_{j=1}^{k-1} x_{i,k-j}^-\otimes \psi^+_{i, j} \ \ (k>  0),\\
\Delta (x_{i,-k}^-)& = x_{i,-k}^-\otimes K_i^{-1}  + 1\otimes x_{i,-k}^- 
+\sum_{j=1}^k  x_{i,-k+j}^-\otimes \psi^-_{i, -j} \ \ (k\ge  0).\end{align*}
\item[(iii)] Modulo $\bu_qX^-\otimes \bu_qX^+$, we have
\begin{equation*} \Delta(h_{i,k}) =h_{i,k}\otimes 1+1\otimes h_{i,k} \ \ 
(k\in\bz).\end{equation*}
\end{enumerate}
\end{prop}
\begin{pf} Part (iii) was proved in \cite{Da}. The rest of the proposition was  
proved in \cite{CPminaff}.
\end{pf}

We conclude this section with some results on the classification of irreducible 
finite--dimensional representations of quantum affine algebras. 
Let \begin{equation*} \cal{A} =\{f\in\bc(q)[[u]]: f(0)=0\}.\end{equation*} 

For any $\bu_q$-module $V$ and any $\mu=\sum_i\mu_i\omega_i\in P$, set
\begin{equation*} V_\mu=\{ v\in V: K_i.v =q_i^{\mu_i}v ,\ \ 
\forall \ i\in I\}.\end{equation*}
We say that $V$ is a module of type 1 if
\begin{equation*} V=\bigoplus_{\mu\in P}V_\mu.\end{equation*}
From now on, we shall only be working with $\bu_q$-modules of type 1.
For $i\in I$,  set
\begin{equation*}
h^\pm_i(u)=\sum_{k=1}^\infty \frac{h_{i,\pm k}}{[k]}u^k.\end{equation*}
 \begin{defn} We say that a $\bu_q$--module $V$ is (pseudo) highest weight, with 
highest weight $(\lambda, \boh^\pm)$, where $\lambda=\sum_{i\in 
I}\lambda_i\omega_i$,  
$\boh^\pm=(\boh_1^\pm(u),\cdots,\boh_n^\pm(u))\in\cal{A}^n$,  if there exists 
$0\ne v\in V_\lambda$ such that $V=\bu_q.v$ and 
\begin{equation*} x_{i,k}^+. v= 0,\ \ K_i.v =q^{\lambda_i}v,\ \ 
h^\pm_i(u).v=\boh_i^\pm(u)v,\end{equation*}
for all $i\in I$, $k\in\bz$.
\hfill\qedsymbol\end{defn}
If $V$ is any highest weight module, then in fact $V=\bu_q(<).v$ and so
\begin{equation*} V_\mu\ne 0\implies \mu=\lambda-\eta\ \ (\eta\in Q^+).
\end{equation*} 
Clearly any highest weight module has  a unique irreducible quotient 
$V(\lambda,\boh^\pm)$.  

The following was proved in \cite{CP}.
\begin{thm}\label{classify} Assume that the pair $(\lambda, \boh^\pm)\in P\times \cal{A}^n$ satisfies, the following: $\lambda=\sum_{i\in I}\lambda_i\omega_i\in 
P^+$, and there exist elements $a_{i,r}\in\bc(q)$ ($1\le r\le \lambda_i, i\in 
I$) such that
\begin{equation*}
\boh_i^\pm(u) =-\sum_{r=1}^{\lambda_i}\text{ln}(1-a_{i,r}^{\pm 1}u). 
\end{equation*}
Then, $V(\lambda,\boh^\pm)$ is the unique (up to isomorphism)  irreducible  
finite--dimensional $\bu_q$--module with highest weight $(\lambda,\boh^\pm)$.
\hfill\qedsymbol\end{thm}
\begin{rem} This statement is actually a reformulation of the statement in 
\cite{CP}.  Setting $\pi_i(u) =\prod_{r=1}^{\lambda_i}(1-a_{i,r}u)$ and 
calculating the eigenvalues of the $\psi_{i,k}$ gives the result stated in 
\cite{CP}. See also \cite{CPweyl}.
\end{rem}
From now on, we shall only be concerned with the modules $V(\lambda, \boh^\pm) $ 
satisfying  the conditions of Theorem \ref{classify}. In view of the preceding 
remark, it is clear that the isomorphism classes of such modules are indexed by 
an $n$--tuple of polynomials $\bpi=(\pi_1,\cdots, \pi_n)$, which have constant 
term 1, and which are split over $\bc(q)$. We shall denote the corresponding 
module by $V(\bpi)$ and the highest weight  vector by $v_\bpi$, where 
$\lambda=\sum_{i=1}^n\text{deg}\pi_i$.
For all $i\in I$, $k\in\bz $, we  have
\begin{equation}\label{qrel1}x_{i,k}^+.v_{\bpi} =0,\ \ K_i.v_{\bpi} 
=q_i^{{\text{deg}}\pi_i}v_{\bpi},\end{equation}
and 
\begin{equation}\label{qrel2} \frac{h_{i,\pm k}}{[k]_i}. v_{\bpi} = \boh_{i, \pm 
k}v_{\bpi},\ \ 
(x_{i,k}^-)^{{\text{deg}}\pi_i+1}.v_{\bpi} =0,\end{equation}
where the $\boh_{i,\pm k}$ are determined from the functional equation
\begin{equation*}{\text {exp}}\left(-\sum_{k\ge 0}{\boh_{i,\pm 
k}}u^k\right) = \pi_i^\pm(u),\end{equation*}
with $\pi_i^+(u)=\pi(u)$ and $\pi_i^-(u) = u^{\text{deg}\pi_i}\pi_i(u^{-1})/
\left.\left(u^{\text{deg}\pi_i}\pi_i(u^{-1}\right)\right|_{u=0}$.

For any $\bu_q$--module $V$, let $V^*$ denote its left dual. Let $- : I\to I$ be the unique  diagram automorphism  such that the irreducible $\frak g$ module  $V(\omega_i)\cong V(\omega_{\overline i})$.    There exists an integer $c\in\bz$ depending only on $\frak g$ such 
that \begin{equation*} V(\bpi)^*\cong V(\bpi^*), \ \ \bpi^* 
=(\pi_{\overline{1}}(q^cu),\cdots ,\pi_{\overline{n}}(q^cu)).\end{equation*}
Analogous statements hold for right duals \cite{CPminaff}.
Recall also, that if a module and its dual are highest weight then they must be 
irreducible.  

Finally,
let $\omega:\bu_q\to\bu_q$ be the algebra automorphism and coalgebra 
anti-automorphism obtained by extending
the assignment $\omega(x_{i,k}^\pm)=-x_{i,-k}^\mp$. If $V$ is any 
$\bu_q$-module, let $V^\omega$ be the pull--back of $V$ through $\omega$. Then,
$(V\otimes V')^\omega\cong (V')^\omega\otimes V^\omega$ and \begin{equation*} 
V(\bpi)^\omega= V(\bpi^\omega)\end{equation*} 
where 
\begin{equation*} \bpi^\omega = (\pi^-_{\overline{1}}(q^\kappa u),\cdots 
,\pi^-_{\overline{n}}(q^\kappa u)),\end{equation*}
for a  fixed integer $\kappa $ depending only on $\frak g$.

We conclude this section with some results  in the case when $\frak g=sl_2$.   

\begin{thm} \label{sl2} \hfill
\begin{enumerate}
\item[(i)]  For $a\in\bc(q)$, $m\in\bz^+$, set  
\begin{equation*}\pi_{m,a}(u)=\prod_{r=1}^m (1-aq^{m-2r+1}u). \ \ 
\end{equation*} The irreducible module $V(m,a)$ with highest weight $\pi_{m,a}$ 
is of dimension $m$ and is irreducible as a $\bu_q^{fin}$--module.
\item[(ii)] For $1\le r\le \ell$, let $a_r\in\bc(q)$ and $m_r\in\bz_+$ be such 
that
\begin{equation*} 
r<s\ \implies\ \frac{a_r}{a_s}\ne q^{m_r-m_s-2p}\ \ (1\le p\le m_r).
\end{equation*} 
The tensor product $V(\bpi_{m_1, a_1})\otimes\cdots\otimes V(\bpi_{m_\ell, 
a_\ell})$ is a  highest weight module with highest weight $\pi_{m_1,a_1}\cdots 
\pi_{m_\ell, a_\ell}$ and   highest weight vector $v_{\pi_{m_1,a_1}}\otimes 
\cdots \otimes v_{\pi_{m_\ell,a_\ell}}$.
\item[(iii)] Assume that $a_1,\cdots , a_\ell\in\bc(q) $ are such that if $r<s$ 
then $a_r/a_s\ne q^{-2}$. The module $W(\pi) = V(1,a_1)\otimes\cdots \otimes 
V(1,a_\ell)$ is the universal finite--dimensional highest weight module with 
highest weight 
\begin{equation*} \pi(u)=\prod_{r=1}^\ell (1-a_ru),\end{equation*}i.e. any other 
finite--dimensional highest weight module with highest weight $\pi$ is a 
quotient of $W(\pi)$.
\item[(iv)] Assume that $\pi=\prod_{j=1}^m(1-b_ju)$ ($b_j\in\bc(q)$) and that 
$V(m_1,a_1)\otimes\cdots\otimes V(m_k,a_k)$ is highest weight. Then, the module
$W(\pi)\otimes V(m_1,a_1)\otimes\cdots\otimes V(m_k,a_k)$ is also highest weight 
if 
\begin{equation*} 
b_s/a_r\ne q^{-1-s},
\end{equation*}
for all $1\le s\le m$ and $1\le r\le k$.

\end{enumerate}
\end{thm}
\begin{pf} Part (i) is  proved  in \cite{CPqa}. Part (ii) is proved as in  Lemma 
4.9 in \cite{CPqa}. In fact, the proof given there establishes the stronger 
result stated here. Part (iii) was proved in \cite{CPweyl}
in the case when $\pi(u)\in\bc[q,q^{-1},u]$. In the general case, choose 
$v\in\bc[q,q^{-1}]$ so that $\tilde{\pi}(u)= \pi(uv)$ has all its roots in 
$\bc[q,q^{-1}]$. Let $\tau_v:\bu_q\to \bu_q$ be the algebra and coalgebra 
automorphism defined by sending $x_k^\pm\to v^kx_k^\pm$.  The pull back of 
$V(\pi)$   through $\tau_v$ is $V(\tilde{\pi})$  and  hence $W(\pi(u))\cong 
W(\pi(uv)$. This proves (iii). Part (iv) is now  immediate.

\end{pf}

Throughout this paper, we shall only work with polynomials in $\bc(q)[u]$  which 
are split and have constant term 1. Let $\pi_{m,a}\in\bc(q)$ be the polynomial 
defined in Theorem \ref{sl2}(i).  It is a simple combinatorial fact [CP1] that 
any such polynomial can be written uniquely as a product 
\begin{equation*}\pi(u)=\prod_{j=1}^s\pi_{m_j, a_j},\end{equation*}
where $m_j\in\bz_+ $, $a_j\in\bc(q)$ and 
\begin{equation*}j< \ell \implies \frac{a_{j}}{a_\ell }\ne  q^{\pm(m_j+m_{\ell 
}-2p)},\ \ 0\le p<\text{min} (m_j,m_\ell).\end{equation*}
If $\pi$ and $\pi'$ are two such  polynomials, then  we say that $\pi>\pi'$ 
if for all $1\le j\le s$, $1\le k\le s'$, we have \begin{equation*} 
\frac{a_j}{a'_k}\ne q^{m_j-m_k'-2p}, \ \ 1\le p\le m_j.\end{equation*} 
This is equivalent to saying that if $a$ is any root of $\pi$ and $1\le k'\le 
s'$, then
\begin{equation*}\frac{a}{a_k'}\ne q^{-1-m_k'}.\end{equation*}
Part (iv) of Theorem \ref{sl2} then says that if $\pi>\pi'$ the module 
$W(\pi)\otimes V(\pi')$ is highest weight.

\section{Braid group action} Let $W$ be the Weyl group of $\frak g$ and let 
$\cal{B}$ be the corresponding braid group. Thus, $\cal{B}$ is the group 
generated by elements $T_i$ ($i\in I$) with defining relations:
\begin{align*} T_iT_j &=T_jT_i,\ \ \text{if}\ \ a_{ij} =0,\\
T_iT_jT_i& =T_jT_iT_j\ \ \text{if}\ \ a_{ij}a_{ji} =1,\\
(T_iT_j)^2&= (T_jT_i)^2,\ \ \text{if}\ \ a_{ij}a_{ji}=2,\\
(T_iT_j)^3&= (T_jT_i)^3,\ \ \text{if}\  \ a_{ij}a_{ji} =3,\end{align*}
where $i,j\in\{1,2,\cdots ,n\}$ and $A=(a_{ij})$ is   the Cartan matrix of 
$\frak g$. 

A straightforward calculation gives the following proposition.
\begin{prop} For all $r\ge 1$, the formulas
\begin{equation*} T_ie_j =e_j-q_i^r\frac{[ra_{ji}]_j}{[r]_j}e_i,\end{equation*} 
define a representation  $\eta_r:\cal{B}\to  \text{end}(V_r)$, where $V_r\cong 
\bc(q)^n$ and $\{e_1,\cdots ,e_n\}$ is the standard basis of $V_r$. Further, 
identifying
\begin{equation*} \cal{A}^n\cong\prod_{r=1}^\infty V_r,\end{equation*}
we get  a representation of $\cal{B}$ on $\cal{A}^n$ given by
\begin{align*}\label{braid} 
(T_i\boh)_j & =\boh_j(u),\ \  {\text{if}}\ a_{ji}=0,\\
(T_i\boh)_j & =\boh_j(u) +\boh_i(qu),\ \  {\text{if}}\ a_{ji}=-1,\\
(T_i\boh)_j & =\boh_j(u) +\boh_i(q^3u) +\boh_i(qu),\ \  {\text{if}}\ 
a_{ji}=-2,\\
(T_i\boh)_j & =\boh_j(u) +\boh_i(q^5u) +\boh_i(q^3u)+\boh_i(qu),\ \  
{\text{if}}\ a_{ji}=-3,\\
 (T_i\boh)_i& =-\boh_i(q_i^2u),\end{align*}
for all $i,j\in I$, $\boh\in\cal{A}^n$.
\hfill\qedsymbol\end{prop}

Let $s_i$, $i\in I$ be a set of simple reflections in $W$. For any $w\in W$, let 
$\ell(w)$ be the length of a reduced expression for $w$.  If 
$w=s_{i_1}s_{i_2}\cdots s_{i_k}$ is a reduced expression for $w$, let 
$T_w=T_{i_1}\cdots T_{i_k}$. It is well--known that $T_w$ is independent of the 
choice of the reduced expression. Given $\boh\in\cal{A}^n$ and $w\in W$, we have
\begin{equation*} T_w\boh =T_{i_1}T_{i_2}\cdots T_{i_k}\boh =\left((T_w\boh)_1, 
\cdots ,(T_w\boh)_n\right).\end{equation*}
We can now prove: 
\begin{prop}\label{twh} Suppose that $w\in W$ and $i\in I$ is such that 
$\ell(s_iw)=\ell(w)+1$. There exists an integer $N\equiv N(i,w,\boh)\ge 0$
and non--negative integers $p_{r,j}$ ($j\in I, 1\le r\le N$) such that
 \begin{equation} (T_w\boh)_i=  \sum_{j\in 
I}\sum_{r=1}^{N}p_{r,j}\boh_j(q^ru),\end{equation}
\end{prop}
\begin{pf} Proceed by induction on $\ell(w)$; the induction clearly starts  at 
$\ell(w)=0$. 
Assume that the result is true for $\ell(w)<k$. If $\ell(w)=k$, write $w=s_jw'$ 
with $\ell(w')=k-1$. Notice that  $j\ne i$ since $\ell(s_iw)=\ell(w)+1$.
We get
\begin{equation*}(T_w\boh)_i= (T_jT_{w'}\boh)_i 
=(T_{w'}\boh)_i(u)+\sum_{s=0}^{|a_{ij}|-1}(T_{w'}\boh)_j(q^{2|a_{ij}|-2s-1}u).
\end{equation*}
If $\ell(s_iw')=\ell(w')+1$, the result follows by induction. If 
$\ell(s_iw')=\ell(w')-1$, we have $w=s_js_iw''$. 
Suppose that  $a_{ij}a_{ji}=-1$. Then,  $\ell(s_jw'')=\ell(w'')+1$ and we get
\begin{equation*} (T_jT_iT_{w''}\boh)_i= (T_iT_{w''}\boh)_i+ 
(T_iT_{w''}\boh)_j(qu) = (T_{w''}\boh)_j(q_iu).\end{equation*}
The result again follows by induction.
The  cases when $a_{ij}a_{ji}=2, 3$ are  proved similarly.  We omit the details.
 \end{pf}

 \section{The main theorem}  Our goal in this section is to obtain a sufficient 
condition for a tensor product of two highest weight representations to be 
highest weight.

Let $V$ be any highest weight finite--dimensional $\bu_q$--module with highest 
weight $\bpi$ (or $(\lambda,\boh^\pm)$ as in Theorem \ref{classify}). For all 
$w\in W$, we have
\begin{equation*} \text{dim}\ V_{w\lambda} = 1.\end{equation*}
If  $s_{i_1}\cdots s_{i_k}$ is a reduced expression for $w$, and $\lambda=\sum_i 
\lambda_i\omega_i$,  set $m_k =\lambda_k$ and define non--negative integers 
$m_j$ (depending on $w$), for $1\le j\le k$, by 
\begin{equation*} s_{i_{j+1}}s_{i_{j+2}}\cdots s_{i_k}\lambda 
=m_j\omega_j+\sum_{i\ne j}m'_i\omega_i.\end{equation*}
Let $v_\lambda$ be the highest weight vector in $V$. For $w\in W$, set 
\begin{equation*}v_{w\lambda} = (x_{i_1,0}^-)^{m_1}\cdots 
(x_{i_k,0}^-)^{m_k}.v_\lambda.\end{equation*}
 If $i\in I$ is  such that $\ell(s_iw)=\ell(w)+1$, then 
\begin{equation}\label{extr1} x_{i,k}^+. v_{w\lambda} =0,\ \ \forall k\in\bz. 
\end{equation} To see this, observe that $w\lambda+\alpha_i$ is not a weight of 
$V$, since $w^{-1}\alpha_i\in R^+$ if $\ell(s_iw)=\ell(w)+1$. It is now easy to 
see that  $v_{w\lambda}\ne 0$,  $V_{w\lambda} =\bc(q)v_{w\lambda}$ and 
\begin{equation*} V=\bu_q.v_{w\lambda}.\end{equation*}
 Since $h_{i,k}V_{w\lambda}\subset V_{w\lambda}$ for all $i\in I$, $k\in\bz$ it 
follows that
\begin{equation*} \frac{h_{i,k}}{[k]}.v_{w\lambda} =\boh_{i,k}^wv_{w\lambda}\ \ 
\forall \ i\in I,\  0\ne k\in \bz,\end{equation*}
where $\boh_{i,k}^w\in\bc(q)$. 
Set
\begin{equation*} 
\boh^w_i(u) = \sum_{k=1}^\infty \boh^w_{i,k}u^k, \ \ \boh^w= (\boh_1^w(u), 
\cdots \boh_n^w(u)). 
\end{equation*} 
Recall that $\boh^1=-(\text{ln}\pi_1(u),\cdots ,\text{ln}\pi_n(u))$.
\begin{prop} \label{extr} If $w\in W$, then 
\begin{equation*} \boh^w = T_w\boh^1.
\end{equation*} 
\end{prop}
\begin{pf} 
We proceed by induction on $\ell(w)$. If $\ell(w)=0$ then $w=id$ and  the result 
follows by definition. Suppose that $\ell(w)=1$, say $w=s_j$. Writing 
$\lambda=\sum_j \lambda_j\omega_j$, we have $v_{s_j\lambda} 
=(x_{j,0}^-)^{\lambda_j}.v_\lambda$. We first show that
\begin{equation}\label{lw1} h_{j,k}(u). v_{s_j\lambda} = - 
\boh_j(q_j^2u)v_{s_j\lambda}  = (T_j\boh)_j(u)v_{s_j\lambda} .\end{equation}
The subspace spanned by the elements $\{(x_{j,l}^-)^r.v_\lambda:0\le r\le 
\lambda_j\}$ is a highest weight module for $\bu_j$, hence 
it is enough to prove \eqref{lw1}   for highest weight representations of 
quantum affine $sl_2$.  In fact it is enough to prove it for the module $W(\pi)$ 
of Theorem \ref{sl2}. Using Proposition \ref{comultip}, we see  that the  
eigenvalue of $h_{i,k}$ on the tensor product of the lowest (and the highest) 
weight vectors is 
just the sum of the eigenvalues values in  each representation. This reduces us 
to the case of  the two-dimensional representation, which is trivial.

Next consider the case $\ell(w)=s_i$, with  $i\ne j$. Recall that \begin{equation*}
[h_{i,r}, x_{j,0}^-] = -\frac{[ra_{ij}]_i}{r} x_{j,r}^- ,\ \ \ \ [h_{j,r}, 
x_{j,0}^-] =-\frac{[2r]_j}{r}x_{j,r}^-.\end{equation*}
Hence,
\begin{align*} h_{i,r}(x_{j,0}^-)^{\lambda_j} &= 
(x_{j,0}^-)^{\lambda_j}h_{i,r}+[ h_{i,r}, (x_{j,0}^-)^{\lambda_j}]\\
&= (x_{j,0}^-)^{\lambda_j}h_{i,r} +\frac{[ra_{ij}]_i}{[2r]_j}[h_{j,r}, 
(x_{j,0}^-)^{\lambda_j}].\end{align*}
This gives 
\begin{align*}\frac{h_{i,r}}{[r]_i}. (x_{j,0}^-)^{\lambda_j}.v_\lambda &= 
\boh_{i,r}(x_{j,0}^-)^{\lambda_j}.v_\lambda 
+\frac{[ra_{ij}]_i}{(q_j^r+q_j^{-r})[r]_i}\left[\frac{h_{j,r}}{[r]_j}, 
(x_{i,0}^-)^{\lambda(h_i)}\right].v_\lambda,\\
&= \boh_{i,r}(x_{j,0}^-)^{\lambda_j}.v_\lambda 
+\frac{[ra_{ij}]_{i}}{(q_j^r+q_j^{-r})[r]_i}(-\boh_{j,r}(x_{j,0}^-)^{\lambda_j}.
v_\lambda +\frac {h_{j,r}}{[r]_j}.\left 
(x_{j,0}^-)^{\lambda_j}\right).v_\lambda,\\&=  
\boh_{i,r} (x_{j,0}^-)^{\lambda_j}.v_\lambda 
-\frac{[ra_{ij}]_{i}}{(q_j^r+q_j^{-r})[r]_i}(\boh_{j,r}
+ q_j^{2r}\boh_{j,r}). (x_{j,0}^-)^{\lambda_j}).v_\lambda,\\
&=(\boh_{i,r}-q_j^r\frac{[ra_{ij}]_i}{[r]_i}\boh_{j,r})(x_{j,0}^-)^{\lambda_j}.v
_\lambda)\end{align*}

This proves the result when $\ell(w)=1$. Proceeding by induction on $\ell(w)$, 
write $w=s_jw'$ with $\ell(w')=\ell(w)-1$. Since
$v_{w\lambda}=(x_{j,0}^-)^{m_j}v_{w'\lambda}$ for some $m_j\ge 0$, the inductive 
step is proved exactly  as in the  case $\ell(w)=1$, with $v_\lambda$ being 
replaced by $v_{w'\lambda}$.
 This completes the proof of the proposition.
\end{pf}

\begin{lem}\label{hwsub}
 Let $V$, $V'$ be finite--dimensional highest weight representations with 
highest weights $\bpi$ and $\bpi'$  and highest weight vectors $v_{\lambda}$ and 
$v_{\lambda'}$ respectively. Assume that $v_{w_0\lambda}\otimes 
v_{\lambda'}\in\bu_q(v_\lambda\otimes v_{\lambda'})$. Then, $V\otimes V'$  is 
highest weight with highest weight vector $v_\lambda\otimes v_{\lambda'}$ and 
highest weight $\bpi\bpi'=(\pi_1\pi_1',\cdots ,\pi_n\pi_n')$.\end{lem}
\begin{pf} It is clear from Proposition \ref{comultip} that the element 
$v_\lambda\otimes v_{\lambda'}$ is a highest weight vector with highest weight 
$\bpi\bpi'$.  It suffices to prove that 
\begin{equation*} V\otimes V'=\bu_q(v_\lambda\otimes 
v_{\lambda'}).\end{equation*}
Since $x_{i,k}^-.v_{w_0\lambda} =0$ for all $i\in I$ and $k\in\bz$, it follows 
from Proposition \ref{comultip} that 
\begin{equation*}\Delta(x_{i,k}^-).(v_{w_0\lambda}\otimes v_{\lambda'}) 
=v_{w_0\lambda}\otimes x_{i,k}^-. v_{\lambda'}.\end{equation*}
Repeating this argument we see that  that $v_{w_0\lambda}\otimes 
V'\subset\bu_q(v_\lambda\otimes v_{\lambda'})$. Now applying  the generators 
$E_{\alpha_i}$, $F_{\alpha_i}$ ($i\in\hat{I}$) repeatedly, we see that 
 $V\otimes V'\subset \bu_q(v_\lambda\otimes v_{\lambda'})$. This proves the 
lemma.
\end{pf}

\begin{lem}\label{slextr} Let $w\in W$ and assume that  $i\in I$ is such that 
$\ell(s_iw)=\ell(w)+1$. Then, $v_{w\lambda}\otimes v_{\lambda'}$ generates a 
$\bu_i$--highest weight module with highest weight $(T_w\boh)_i\boh_i'$. 
\end{lem}

\begin{pf} This is immediate from Proposition \ref{comultip} and Proposition 
\ref{extr}. \end{pf}
 We can now prove our main result. Given $\bpi=(\pi_1,\cdots ,\pi_n)$, and $w\in 
W$, set 
\begin{equation*}T_w\bpi=(\text{exp}-(T_w\ln\pi_1(u))_1, \cdots 
,\text{exp}-(T_w\ln\pi_n(u))_n).\end{equation*}

\begin{thm} The module $V(\bpi_1)\otimes \cdots\otimes V(\bpi_r)$ is highest 
weight if for all $w\in W$ and $i\in I$  with $\ell(s_iw)=\ell(w)+1$, we have 
\begin{equation*}  j<\ell \implies\ (T_w\bpi_j)_i> (\bpi_\ell)_i. 
\end{equation*}
\end{thm}
\begin{pf} First observe that by Proposition \ref{twh},  $(T_w\bpi_j)_i$ is 
indeed a polynomial.
If $\frak g=sl_2$, then this is the statement of Theorem \ref{sl2} (ii).  
For arbitrary $\frak g$, proceed by induction on $r$. If $r=1$ there is nothing 
to prove. 
Let $r>1$ and let $V'=V(\bpi_2)\otimes \cdots\otimes V(\bpi_r)$. Then $V'$ is 
highest weight module with highest weight vector 
$v'=v_{\bpi_2}\otimes\cdots\otimes v_{\bpi_r}$ and highest weight 
$\bpi'=\bpi_2\cdots \bpi_r$. Setting $\lambda= (\text{deg}\pi_i,\cdots 
,\text{deg}\pi_n)$, it is enough by Lemma \ref{hwsub} to prove that 
\begin{equation*} v_{w_0\lambda}\otimes v'\in \bu_q(v_{\bpi_1}\otimes 
v').\end{equation*}
Writing $w_0=s_{i_1}s_{i_2}\cdots s_{i_N}$ and using  Lemma \ref{slextr}, it 
suffices to prove that for all $1\le j\le N$, 
\begin{equation*} v_{s_{i_j}s_{i_{j+1}}\cdots s_{i_N}\lambda}\otimes 
v'\in\bu_{i_j}(v_{s_{i_{j+1}\cdots s_{i_N}}\lambda}\otimes v').\end{equation*}
or equivalently
that the $\bu_{i_j}$--module $\bu_{i_j}.v{s_{i_{j+1}\cdots 
s_{i_N}}\lambda}\otimes \bu_{i_j}.v'$ is highest weight. Taking $w=id$, we have 
$(\bpi_l)_{i_j}>(\bpi_s)_{i_j}$ if $l<s$. It thus  follows from Theorem 
\ref{sl2}(ii)  that 
\begin{equation*} 
\bu_{i_j}.v'=\bu_{i_j}.v_{\bpi_2}\otimes\cdots\otimes\bu_{i_j}.v_{\bpi_r}.
\end{equation*}
Since $(T_w\bpi_1)_i> (\bpi_\ell)_i$, the result follows from 
 Theorem \ref{sl2}(iv).
\end{pf}

\section{Relationship wth Kashiwara's results and conjectures}

Let us consider the special case when $\boh$ has the following form,
\begin{equation*} \boh_j^\pm(u)=0, \ \ j\ne i, \ \ \boh_i^\pm (u) 
=-\sum_{r=1}^m\text{ln}(1-aq_i^{m-2r+1}u),\end{equation*}
and denote the corresponding $n$--tuple of power series  by $\boh^i_{m,a}$ and 
the $n$--tuple of polynomials by $\bpi^i_{m,a}$. 
We shall prove the following result.

\begin{thm}\label{k} Let $i_1,i_2\cdots, i_l\in I$, $a_1,\cdots, a_l\in\bc(q)$, 
$m_1,\cdots, m_\ell \in\bz_+$, and assume that
\begin{equation*} r<s\ \ \implies \ \ \frac{a_r}{a_s}\ne 
q^{d_{i_r}m_r-d_{i_s}m_s- d_{i_r}-d_{i_s}- p}\ \ \forall\ \ p\ge 0. 
\end{equation*}
Then, the tensor product $V(\bpi^{i_1}_{m_1,a_1})\otimes\cdots\otimes 
V(\bpi^{i_\ell}_{m_\ell,a_\ell})$ is a highest weight module.
\end{thm}
Assume the theorem for the moment.

\begin{rem}
In the special case when $m_{j}=1$ for all $j$, it was conjectured in \cite{AK} 
that such  a tensor product is cyclic if $a_{j}/a_{l}$ does not have a pole at 
$q=0$ if $j<\ell$, and this was proved when $\frak g$ is of type $A_n$ or $C_n$; 
subsequently,  a geometric proof of this conjecture  was given in \cite{VV} when 
$\frak g$ is simply--laced; a complete proof was given using crystal basis 
methods in \cite{K}. \end{rem}

The following corollary to Theorem  \ref{k} was conjectured in \cite{K}. 
\cite{HKOTY}.

\begin{cor}
 The  tensor product $V=V(\bpi^{i_1}_{m_1, 1})\otimes\cdots\otimes 
V(\bpi^{i_\ell}_{m_\ell,1})$ is an irreducible $\bu_q$--module.
\end{cor}
\begin{pf} First observe that if $d_{i_1}m_1\le d_{i_2}m_2\le \cdots \le 
d_{i_l}m_\ell$ then the tensor product is cyclic by Theorem \ref{k}. We claim 
that it suffices to prove the corollary in the case when $\ell=2$.  For then,  
by rearranging the factors in the tensor product we can show that both $V$ and 
its dual are  highest weight and hence  irreducible.  To see that  
$V=V(\bpi^{i_1}_{m_1, 1})\otimes V(\bpi^{i_2}_{m_2, 2})$ is cyclic if 
$d_{i_1}m_1> d_{i_2}m_2$, it is enough to show that  $V^\omega$ is cyclic, since 
$\omega$ is an algebra  automorphism.
Now,  $V^\omega = V(\bpi^{i_1}_{m_2, v^{-1}})\otimes V(\bpi^{i_1}_{m_1, 
v^{-1}})$
for some fixed $v$ depending only on $\frak g$, and this is cyclic by the 
theorem. This proves the result.

 \end{pf}

It remains to prove the theorem, for which we must show that, if $j<\ell$  and 
$\ell(s_iw)=\ell(w)+1$, then
\begin{equation} \label{crux}\left(T_w\bpi^{i_j}_{m_j,a_j}\right)_i >  
(\bpi^{i_\ell}_{m_\ell,a_\ell})_i.\end{equation}

Using Propositon \ref{twh}, we see that
\begin{equation*} \left(T_w\bpi^{i_j}_{m_j,a_j}\right)_i =\prod_{r\ge 
0}\pi_{m_j,1}(q^ru),\end{equation*}
where $r$ varies over a finite subsubset of $\bz_+$ with multiplicity. 
This means that  any root of $(T_w\bpi^{i_j}_{m_j,a_j})_i$ has the form 
$q^{d_{i_j}(m_j-2p+1)+r}a_j$ where $r\ge 0$ and $1\le p\le m_j$, and hence, 
using the assumption on $a_j/a_\ell$, that
\begin{equation*}  \frac{q^{d_{i_j}(m_j-2p+1)+r}a_j}{a_l}\ne 
q^{-1-m_\ell}.\end{equation*}
This   proves \eqref{crux} and the proof of the theorem is complete.

\end{document}